\documentclass{amsart}

\usepackage{graphicx}
\usepackage[export]{adjustbox}
\usepackage{subfigure}
\usepackage{caption}

\usepackage{textcomp}

\usepackage{amsfonts, amsmath, amssymb,amsthm,verbatim}  
\usepackage{times,enumerate}
\usepackage{nccmath}
		{\end{pmatrix}\end{medsize}}%
\newenvironment{mbmatrix}{\begin{medsize}\begin{bmatrix}}%
		{\end{bmatrix}\end{medsize}}%
\usepackage{rotating}
\usepackage{lscape}
\usepackage{multicol}

\usepackage{thmtools, thm-restate}
\usepackage{blkarray}
\usepackage{hyperref}

\usepackage{color}
\DeclareGraphicsExtensions{.pdf,.png,.jpg}

\usepackage{algorithmic}
\usepackage{algorithm2e}
\RestyleAlgo{boxruled}

\usepackage{float}

\newcommand{\benu}{\begin{enumerate}}
		\newcommand{\eenu}{\end{enumerate}}
\newcommand{\beqn}{\begin{align*}}
	\newcommand{\eeqn}{\end{align*}}
\newcommand{\bdefi}{\begin{definition}}
	\newcommand{\edefi}{\end{definition}}
\newcommand{\bcor}{\begin{corollary}}
	\newcommand{\ecor}{\end{corollary}}
\newcommand{\bthe}{\begin{theorem}}
	\newcommand{\ethe}{\end{theorem}}
\newcommand{\bpro}{\begin{proposition}}
	\newcommand{\epro}{\end{proposition}}
\newcommand{\blem}{\begin{lemma}}
	\newcommand{\elem}{\end{lemma}}
\newcommand{\brem}{\begin{remark}}
	\newcommand{\erem}{\end{remark}}
\newcommand{\bequ}{\begin{equation}}
    \newcommand{\eequ}{\end{equation}}
\newcommand{\bprf}{\begin{proof}}
	\newcommand{\eprf}{\end{proof}}

\newtheorem{theorem}{Theorem}[section]
\newtheorem{corollary}[theorem]{Corollary}
\newtheorem{lemma}[theorem]{Lemma}
\newtheorem{proposition}[theorem]{Proposition}
\newtheorem{criterion}[theorem]{Criterion}

\theoremstyle{definition}
\newtheorem{definition}[theorem]{Definition}
\newtheorem{example}{Example}

\newcommand{\mc}{\mathcal}
\newcommand{\mbb}{\mathbb}

\newcommand{\norm}[1]{\left\lVert#1\right\rVert}

\theoremstyle{remark}
\newtheorem{remark}[theorem]{Remark}

\numberwithin{equation}{section}

\begin{document}

\title{Practical construction of positive maps which are not completely positive}

\author{Abhishek Bhardwaj}
\address{MATHEMATICAL SCIENCES INSTITUTE, THE AUSTRALIAN NATIONAL UNIVERSITY, UNION LANE, CANBERRA ACT 2601}
\curraddr{}
\email{Abhishek.Bhardwaj@anu.edu.au}
\thanks{}

\subjclass[]{}

\date{}

\dedicatory{}


\begin{abstract}

	This article introduces	PnCP, a MATLAB toolbox for constructing positive maps which are not completely positive. We survey optimization and sum of squares relaxation techniques to find the most numerically efficient methods for this construction. We also show how this package can be applied to the problem of classifying entanglement in quantum states.

\end{abstract}

\maketitle

\section{Introduction}

For $n\in\mbb{N}$, let $M_{n}(\mbb{R})$ be the vector space of $n\times n$ matrices over $\mbb{R}$. A matrix $A\in M_{n}(\mbb{R})$ is called \textit{positive semi-definite} if for all $v\in\mbb{R}^{n}$, $v^{T}Av \geq0$; in this case we write $A\succeq0$. Given two matrix spaces $\mc{A}$ and $\mc{B}$, a linear map $\Phi : \mc{A} \rightarrow \mc{B}$ with the involution-preserving property $\Phi(A^{*}) = \Phi(A)^{*}$ for all $A\in\mc{A}$, is called \textit{positive} if for all $A\succeq0$, $\Phi(A)\succeq0$.
For a given $k\in\mbb{N}$, such linear maps induce the ampliation
$$
\Phi^{(k)} : M_{k}(\mbb{R})\otimes\mc{A} \rightarrow M_{k}(\mbb{R})\otimes\mc{B} ;\qquad M\otimes A \rightarrow M\otimes\Phi(A),
$$
where $\otimes$ is the standard Kronecker tensor product of matrices. If $\Phi^{(k)}$ is positive then we call $\Phi$ \textit{k-positive}. If $\Phi^{(k)}$ is positive for all $k\in\mbb{N}$, then $\Phi$ is called \textit{completely positive}. Positive and completely positive maps arise naturally in matrix theory and operator algebras (e.g., positive linear functionals) \cite{paulsen2002completely, woronowicz1976positive}, frequently in quantum information theory \cite{horodecki2003entanglement, szarek2010often, nielsen2000quantum}, and have recently even been used in semi-definite programming \cite{klep2013exact}. 

We study these maps via their correspondence to positive and non-negative polynomials. Let $\mbb{S}_{n}$ be the subspace of symmetric matrices $\{ A\in M_{n}(\mbb{R}) : A^{T} = A \}$. Restricting these involution-preserving maps to the space of symmetric matrices, each $\Phi :\mbb{S}_{n} \rightarrow \mbb{S}_{m}$ gives rise to a biquadratic, bihomogeneous polynomial $p_{\Phi}\in\mbb{R}[x,y]$, with 
$$
p_{\Phi}(x,y) = y^{T}\Phi(x^{T}x)y.
$$
It is easily seen (or see, e.g.,  \cite{klep2017there}) that $\Phi$ is positive if and only if $p_{\Phi}$ is non-negative on $\mbb{R}^{n+m}$, and $\Phi$ is completely positive if and only if $p_{\Phi}$ is a sum of squares (SOS) on $\mbb{R}^{n+m}$. 

The connection between non-negative and SOS polynomials plays a central role in real algebraic geometry. There are many results concerning this interplay, see for instance the surveys \cite{belton2019panorama, powers2011positive, laurentsurvey1, laurentsurvey2} or the book \cite{marshall2008positive}. In particular, \cite{blekherman2016sums} explores the connection between varieties of minimal degrees and non-negative polynomials. Their main theorem (given below) shows that on varieties of minimal degrees, non-negative quadratic forms have an SOS decomposition with linear forms.
\begin{theorem}[Thereom 1.1, \cite{blekherman2016sums}]\label{BSV}
	Let $X \subseteq \mbb{P}^{n}$ be a real irreducible non-degenerate projective sub-variety, with homogeneous coordinate ring R, such that the set $X(\mbb{R})$ of real points is Zariski dense. Every non-negative real quadratic form on $X$ is a sum of squares of linear forms in $R$ if and only if $X$ is a variety of minimal degree.
\end{theorem}
\noindent Moreover, when $X$ is not of minimal degree, \cite{blekherman2016sums} gave a construction for generic quadratic forms which are non-negative on $X$ but not SOS. In \cite{klep2017there} the authors specialize this construction (Procedure 3.3 of \cite{blekherman2016sums}) to biquadratic, bihomogeneous polynomials (biforms) over the Segre Variety, which is the image of the Segre embedding $\sigma:\mbb{P}^{n}\times\mbb{P}^{m}\rightarrow\mbb{P}^{(n+1)(m+1)-1}$ (and is well known to not be of minimal degree for $n,m\geq3$). This formalization of the method in \cite{blekherman2016sums}, gives an algorithmic construction of positive maps which are not completely positive (pncp maps for short). Letting $n,m>2$, $d=n+m-2$ and $N=n+m$, the algorithm of \cite{klep2017there} (\cite{blekherman2016sums}) can be summarized as:
\begin{algorithm}
	\caption{KMSZ construction}\label{alg1}
	1. Generate random points $x\in\mbb{R}^{n}, y\in\mbb{R}^{m}$\\
	2. Use $x,y$ to create bilinear forms $\{ h_{0}, \dotso, h_{d} \}$ over $\mbb{R}^{N}$ \\
	3. Generate $f \notin \langle h_{0}, \dotso, h_{d} \rangle$ so that $f \neq SOS$ on $\mbb{R}^{N}$\\
	4. Choose $\delta$ small enough so that $F_{\delta} = \delta f + h_{0}^{2}+\dotsb+h_{d}^{2} \geq 0$ on $\mbb{R}^{N}$
\end{algorithm}

\noindent Steps 1-3 are simple linear algebra computations, our contribution in this work is to find the most practical technique for Step 4, and to establish benchmarks for this type of construction.


This is an expository and experimental article in which we introduce the MATLAB package PnCP, currently the only implementation of Algorithm \ref{alg1}. 
We survey recent optimization techniques for verifying Step 4 and specify relaxations theoretically superior to those presented in \cite{klep2017there}. We implement and test these methods in PnCP. Our package and test data are made available at \url{https://bitbucket.org/Abhishek-B/pncp/}. 
We also consider rationalizations of the forms obtained with Algorithm \ref{alg1} to obtain exact certificates of non-negativity (PnCP is able to construct pncp maps with rational coefficients).

PnCP is developed as a consequence of the rising interest in quantum information and its purpose is to help identify entangled (quantum) states; pncp maps preserve their positivity on separable states, however they may fail to preserve positivity on entangled states, which provides the following classification criterion.
\begin{criterion}[The general criteria, \cite{audretsch2008entangled} section 8.4]
	A quantum state $\rho\in M_{t}(\mbb{R})$ is entangled if there is a pncp map $\Phi$ such that the ampliation $(I\otimes \Phi)(\rho)\nsucceq0$.
\end{criterion}

As an example, consider the Bell State, which has density matrix (see Section \ref{quantum})
$$
\rho = \frac{1}{2} \begin{bmatrix}
1 & 0 & 0 & 1 \\
0 & 0 & 0 & 0 \\
0 & 0 & 0 & 0 \\
1 & 0 & 0 & 1 
\end{bmatrix} \in M_{2}(\mbb{R})\otimes M_{2}(\mbb{R}).
$$
and let $\Phi$ be the transposition map (clearly positive, and known to be pncp). Then the ampliation $(I\otimes\Phi):M_{2}(\mbb{R})\otimes M_{2}(\mbb{R})\rightarrow M_{2}(\mbb{R})\otimes M_{2}(\mbb{R})$ applied to $\rho$ gives,
$$
(I\otimes\Phi)(\rho) = \frac{1}{2}\begin{bmatrix}
1 & 0 & 0 & 0 \\
0 & 0 & 1 & 0 \\
0 & 1 & 0 & 0 \\
0 & 0 & 0 & 1 
\end{bmatrix} 
$$
which has a negative eigenvalue of $-1/2$, and serves as evidence of entanglement in the Bell State. While the transposition map was sufficient in this simple example, in general finding a suitable map is difficult. With the help of PnCP one can generate many such maps to test for entanglement (see the examples in Section \ref{quantum} for details).

The article is organized as follows. 
Section 2 reviews some notation and background for the optimization involved in Step 4. 
In Section 3 we present some of the relaxations we surveyed and thought to be promising for using in Step 4. We also present our implementation of these methods using MATLAB and show their performance via computational efficiency (w.r.t.\ time) and success rate.
Section 4 details issues in generating pncp maps with rational coefficients using Algorithm \ref{alg1}. We also show the difference in computational requirements for constructing maps with floating point coefficients and those with rational coefficients. 
Section 5 explains how we use PnCP to identify entanglement in quantum states. We demonstrate this usefulness through illustrative examples. 

\subsection*{Acknowledgement} I wish to thank Prof. Igor Klep and Dr. Alja\v{z} Zalar for introducing me to this topic, and for insightful discussions throughout the project.

\section{Background}
In this section we present the necessary mathematical background and notation for undertaking Step 4 of Algorithm \ref{alg1}. We focus on the general optimization problem of Step 4, and the underlying principles for finding a solution.

We first look at minimization techniques which we can use to ensure non-negativity, and then consider their relaxations which make them computationally feasible. We also describe how we implement these techniques in PnCP for practical success. 

We use the following notation; $\mbb{N}$ (resp. $\mbb{R}, \mbb{C}$) denotes the usual set of non-negative integers (resp. real numbers, complex numbers). We write $\mbb{R}[x] = \mbb{R}[x_{1},\dotso,x_{n}]$ for the ring of real polynomials in $n$ variables and $\mbb{R}[x]_{d}$ for the subset of polynomials in $\mbb{R}[x]$ with total degrees bounded by $d$. We occasionally also write $\mbb{R}[x,y] = \mbb{R}[x_{1},\dotso,x_{n},y_{1},\dotso,y_{m}]$ for the special case of polynomials in two sets of variables. For any integer $n>0$, $[n] = \{1,\dotso,n\}$ and for a subset $I\subseteq [n]$, $|I|$ denotes its cardinality. For $k\in\mbb{N}$, $[n]_{k} = \{ I\subseteq[n] : |I| = k \}$. 

A subset $I\subseteq\mbb{R}[x]$ is called an \textit{ideal} if $\mbb{R}[x] \cdot I \subseteq I$. The set $\langle g_{1}, \dotso, g_{j} \rangle$ is the ideal generated by $\{ g_{1},\dotso,g_{j} \} \subseteq \mbb{R}[x]$, which is the smallest ideal containing $\{ g_{1},\dotso,g_{j} \}$. According to the \textit{Hilbert Basis Theorem} \cite{cox2007ideals}, every ideal has such a finite generating set. The \textit{variety} of an ideal is the set of common complex zeros for the ideals' generators,
$$
V(I) = V(\langle g_{1},\dotso,g_{j}\rangle) = \{ x\in\mbb{C}^{n} \ : \ g_{k}(x)=0,\ \forall k=1,...,j \},
$$
or more generally 
$$
V(I) = \{ x\in\mbb{C}^{n} \ : \ p(x)=0, \ \forall p\in I \}.
$$
The \textit{real variety} of $I$ is simply the restriction of $V(I)$ to the reals. We denote this with $V^{\mbb{R}}(I)$. If the variety $V(I)$ is a finite set, then $I$ is called \textit{zero dimensional} (this is not the same as requiring $V^{\mbb{R}}(I)$ to be finite). For every ideal $I\in\mbb{R}[x]$, its \textit{radical} is the ideal
$$
\sqrt{I} = \{ p\in\mbb{R}[x] \ : \ p^{r}\in I \text{ for some } r\in\mbb{N} \}.
$$
For more details see \cite{cox2007ideals}. With any finite set $G = \{g_{1},\dotso,g_{j} \}\subseteq\mbb{R}[x]$ we have the \textit{semialgebraic} set, and the \textit{preorder} generated by $G$ resp.,
\begin{equation*}
\begin{gathered}
S(G) = \{ x\in\mbb{R}^{n} \ : \ g_{k}(x)\geq0,\ \forall k=1,...,j \}, \\
PO(G) = \left\{ \sum_{\gamma\in\{0,1\}^{j}} s_{\gamma}g_{1}^{\gamma_{1}}\dotsb g_{j}^{\gamma_{j}} \ : \ s_{\gamma} \text{ is SOS in } \mbb{R}[x] \right\}.
\end{gathered}
\end{equation*}

\subsection{Minimization} 
Consider a general constrained minimization problem
\begin{equation}\label{minprob}
	\begin{gathered}
		\min_{x\in\mbb{R}^{N}} p(x) \\
		s.t. \ \ g_{1}(x) = \dotsb g_{t}(x) = 0.
	\end{gathered}
\end{equation}
The difficulty in Step 4 consists of solving a system like this and verifying the solution to be non-negative.  However, as is well known, testing the non-negativity of a polynomial $p$ is an NP-hard problem \cite{murty1987some, parrilo2003semidefinite}. Instead, we generally use an \textit{SOS relaxation} which is computationally tractable (the idea is to decompose $p$ as $p=\sum_{i}q_{i}^{2}$ modulo the constraints). For the generic minimization problem (\ref{minprob}), the standard relaxation is given by 
\begin{equation}\label{sosprog}
	\begin{gathered}
		\max{\gamma}, \\ 
		s.t. \ \ p(x) - \gamma = SOS + \sum_{1}^{t}\phi_{j}(x)g_{j}(x), \\
		\phi_{j}\in\mbb{R}[x],\\
		deg(\phi_{j}g_{j})\leq 2k.
	\end{gathered}
\end{equation}
Letting $p_{min}$ and $\gamma_{k}$ be the solutions to (\ref{minprob}) and (\ref{sosprog}) respectively, Lasserre \cite{lasserre2001global} has shown that $\gamma_{k}\rightarrow p_{min}$ as $k\rightarrow \infty$ under some natural conditions.

The recommended relaxation in \cite{klep2017there} for $F_{\delta}$ in Step 4 of Algorithm \ref{alg1} is \\
\begin{equation} 
	\label{CNR}
	\begin{gathered}
		\max_{\delta>0}{\delta},  \\
		s.t. \  \left(\sum (x_{i} y_{j})^{2}\right)^{\ell}F_{\delta} (x,y) = SOS, \\
		\ell \in \mbb{N}.
	\end{gathered}
\end{equation}
These relaxation problems can be stated and solved as an appropriate optimization program (semi-definite, second order cone, quadratically constrained, etc.). In recent years, there have been many developments in optimization for computing minima, and the majority of solvers can handle the broad class of these problems.

\section{Relaxations \& Performance}\label{relaxations}

We now present alternate SOS relaxations to solving problem (\ref{minprob}). We present the theory in this section with regards to an arbitrary function $p(x)\in\mbb{R}[x]$. We then give a description of how the results apply to our function of interest $F_{\delta}$, and finally discuss the implementation and performance. 

\subsection{Rational Functions} Let us begin by considering Artin's solution to Hilbert's $17^{th}$ problem \cite{bochnak2013real}.
\begin{theorem}[Hilbert's $17^{th}$ Problem]\label{artinsol}
	For any $p\in\mbb{R}[x]$, if $p\geq0$ on $\mbb{R}^{n}$, then $p$ is a sum of squares of rational functions, i.e., there are polynomials $g, q_{i} \in\mbb{R}[x]$, with $g\neq0$, such that $$g^{2}p = \sum_{i}q^{2}_{i}.$$
\end{theorem}
This result provides the most fundamental SOS relaxation. For Step 4, instead of minimizing $F_{\delta}$, we look for a decomposition into sums of rational squares, i.e.,
\begin{equation} 
\label{HilB}
	\begin{gathered}
		\max_{\delta>0}{\delta},  \\
		s.t. \ \sigma(x,y)F_{\delta} (x,y) = SOS, \\
		\sigma(x,y) \in \mbb{R}[x,y],
	\end{gathered}
\end{equation}
with $\sigma(x,y)=SOS$. If for some $\delta$, $F_{\delta}$ is non-negative, then by Theorem \ref{artinsol} the SOS decomposition in (\ref{HilB}) always exists. 

Note that (\ref{HilB}) is a quadratically constrained optimization program (non-linear in the decision variables, $\delta$ and the coefficients of $\sigma$), which can be solved with solvers such as PENLAB \cite{Penlab}, but our early tests indicated that this approach is not ideal. So we instead implement (\ref{HilB}) with a ``bisection" approach. This is already the suggested method in \cite{klep2017there}, which tried to solve (\ref{CNR}), and increases $\ell$ if a solution is not found. While bisecting may not find the optimal $\delta$, it does find a successful instance of $\delta$ very quickly.

For the Hilbert method (\ref{HilB}), let $G$ be the Gram matrix of $\sigma$. We fix $\delta = 2^{0}$, $k=1$, and solve the following
\begin{equation*} 
\begin{gathered}
\text{find } \sigma(x,y) \in \mbb{R}[x,y]_{k},\\
s.t. \ Tr(G)=1, \\
\sigma(x,y)F_{\delta} (x,y) = SOS.
\end{gathered}
\end{equation*}
If a solution is not found, we first bisect over $\delta$, and if still there is no solution we increase $k$ and repeat. We set the limits of $\delta$ to be $2^{-6}$ and $k$ to be $2$.

The SOS decomposition and related optimization problems are generated using the symbolic computation package YALMIP {\cite{Lofberg2004, Lofberg2009}}. Our MATLAB code for the experiments, as well as our computational data is available on \url{https://bitbucket.org/Abhishek-B/pncp/}, so that the reader may verify the results of our experiments. 

To solve the required SDP we use the MOSEK solver \cite{mosek} with our implementations. Verification of the SOS decomposition is done with the YALMIP command \texttt{sol.problem==0} (where \texttt{sol} is what we name our solution), as well as requiring the residual of the problem to be small ($\leq O(10^{-6})$). 

All of the experiments were carried out on a standard Dell Optiplex 9020, with 12GB of memory, an Intel \textregistered Core™ i5-4590 CPU @ 3.30GHz$\times$4 processor, 500GB of storage and running Ubuntu 18.04 LTS.

\graphicspath{{Images/}}
\begin{figure}[H]
	\includegraphics[width=0.9\linewidth, center]{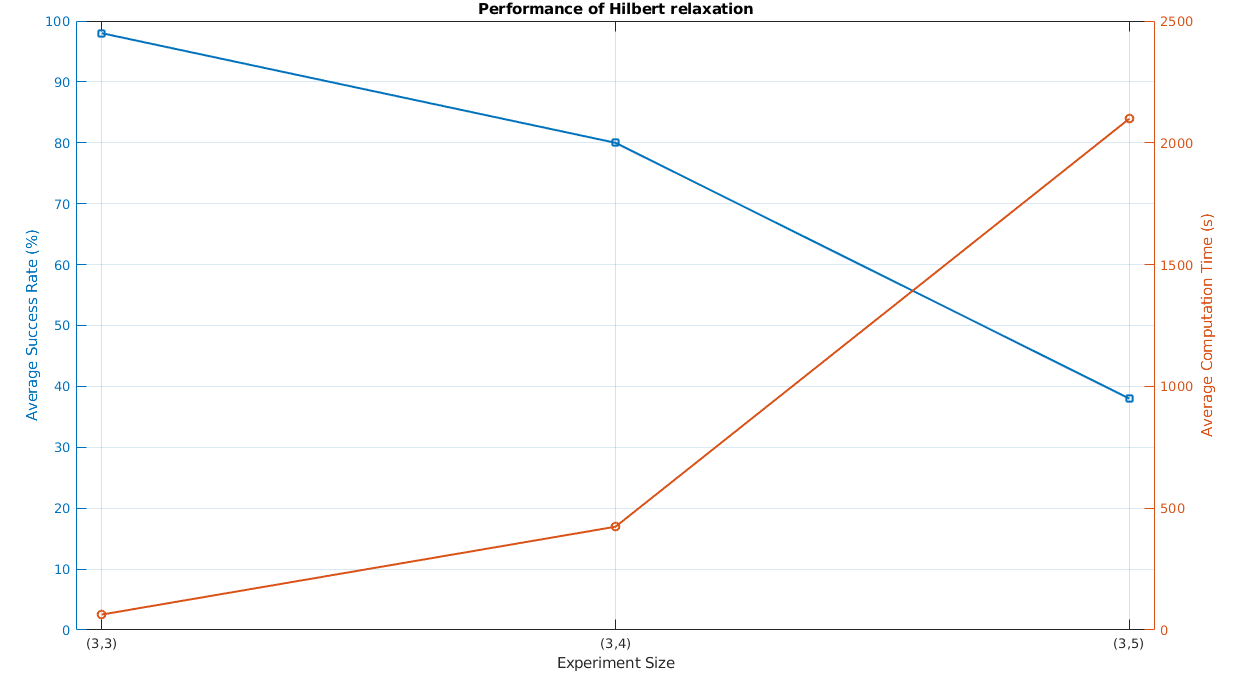}
	\caption{}
	\label{fig1}
\end{figure}

The success rate of this relaxation for problems of small size is remarkable, as seen in Figure \ref{fig1}. Moreover, we observe from the average residual (which includes the failed examples as well) in Table \ref{table:1}, that if we were to allow the residual to be slightly larger (say $\leq O(10^{-5})$), we would see a higher success rate. This would also reduce computation times, increasing the appeal of this relaxation.
\begin{table}[H]
	\centering
	\begin{tabular}{ |p{2cm}||p{2cm}|p{2cm}|p{2cm}|  }
		\hline
		\multicolumn{4}{|c|}{Hilbert Relaxation}    \\
		\hline
		$(n,m)$ & Success (\%)  & Time (s)  & Residual    \\
		\hline
		$(3,3)$ & $98$          & $63.31$   & $7.19\times 10^{-7}$     \\
		$(3,4)$ & $80$          & $423.99$  & $2.02\times 10^{-6}$     \\
		$(3,5)$ & $38$          & $2098.93$ & $1.17\times 10^{-5}$     \\
		\hline
	\end{tabular}
	\vspace{2pt}
	\caption{Average performance of relaxation (\ref{HilB})}
	\label{table:1}
\end{table}

\begin{remark}
	\label{initial}
	After running a few experiments it becomes apparent that in the Hilbert method, we should initialize $k=2$. While there are instances where $k=1$ has a solution, it works with very small $\delta$ and hence requires a long runtime due to the number of bisections.
	We also add $Tr(G)=1$ in our constraints to avoid the trivial solution of $\sigma \equiv 0$. 
\end{remark}

The relaxation (\ref{CNR}) is a simplified version of (\ref{HilB}), which fixes the denominator
\begin{equation*} 
\sigma(x,y) = \left(\sum (x_{i} y_{j})^{2}\right)^{\ell}.
\end{equation*}
We refer to this simplification as the \textit{Coordinate Norm Relaxation} (CNR) and implement it similar to the Hilbert method. Since $\sigma$ is known, we maximize $\delta$ and ``bisect" over $\ell\leq 2$. The verification of a solution is also similar, with the additional requirement $\delta > O(10^{-4})$ as otherwise $\delta$ becomes indistinguishable from numerical error.

\graphicspath{{Images/}}
\begin{figure}[H]
	\includegraphics[width=0.9\linewidth, center]{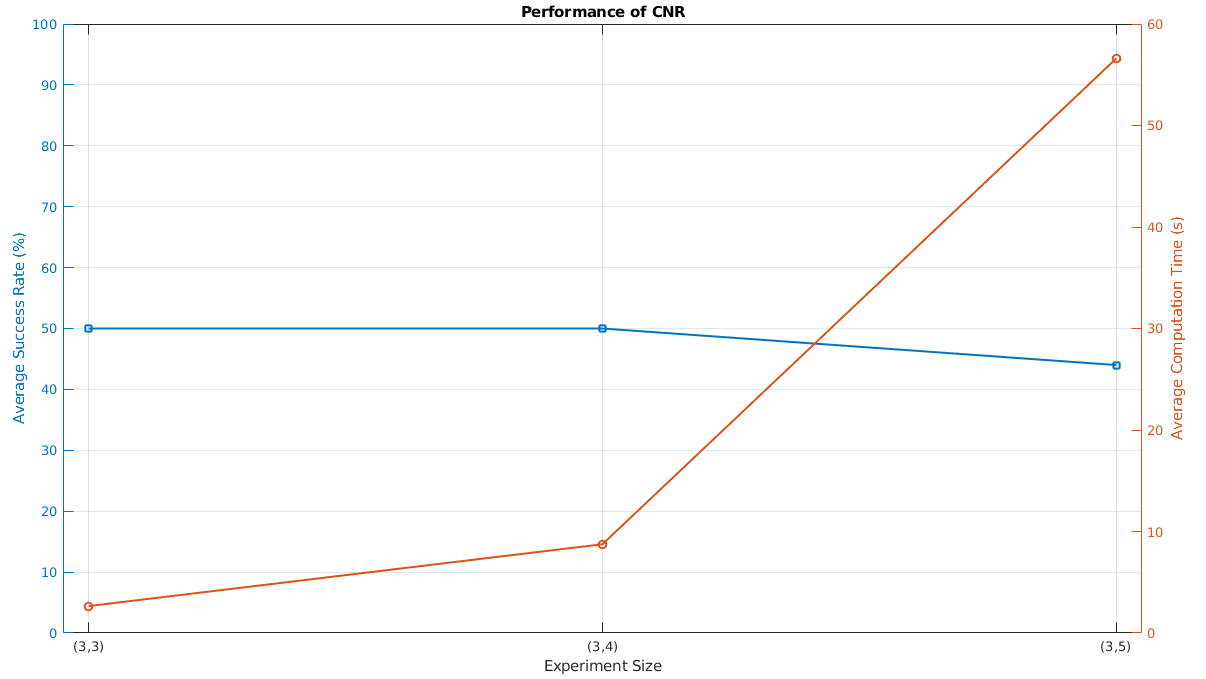}
	\caption{}
	\label{fig2}
\end{figure}

As we can see (Figure \ref{fig2} or Table {\ref{table:2}}), this relaxation is incredibly fast (it is in fact the fastest relaxation). On problems of smaller size, it is not as successful compared to the Hilbert method, but we can see from the residuals, that if we relax our verification criteria, we might improve the success rate of the CNR quite dramatically.
\begin{table}[H]
	\centering
	\begin{tabular}{ |p{2cm}||p{2cm}|p{2cm}|p{2cm}|p{2cm}|  }
		\hline
		\multicolumn{5}{|c|}{CNR}    \\
		\hline
		$(n,m)$ & Success (\%)  & Time (s) & Residual               & Average $\delta$  \\
		\hline
		$(3,3)$ & $50$          & $2.65$   & $4.89\times 10^{-6}$   &  1.83   \\
		$(3,4)$ & $50$          & $8.75$   & $5.53\times 10^{-6}$   &  0.13   \\
		$(3,5)$ & $44$          & $56.61$  & $1.34\times 10^{-5}$   &  0.09   \\
		\hline
	\end{tabular}
	\vspace{2pt}
	\caption{Average performance of relaxation (\ref{CNR})}
	\label{table:2}
\end{table}

If we consider the variables $z_{ij} = x_{i}\otimes y_{j}$ over the Segre variety, then the CNR can be written as
\begin{equation*} 
\label{CNR_z}
\begin{gathered}
\max_{\delta>0}{\delta},  \\
s.t. \  \left(\sum z_{ij}^{2}\right)^{\ell}F_{\delta} (z) = SOS, \\
\ell \in \mbb{N}.
\end{gathered}
\end{equation*}
For strictly positive polynomials $p$, there always exists an $\ell$ such that the denominator $(x_{1}^{2}+\dotsb+x_{n}^{2})^{\ell}$ allows an SOS decomposition (see the second and third Theorems of \cite{reznick1995uniform}). However the appropriate choice of $\ell$ depends on the minimum of $p$, and in fact $\ell \rightarrow \infty$ as $\inf(\{p(x) | x\in\mbb{R}^{n}\})\rightarrow0$. For polynomials with zeros, this denominator has been used in practice (see \cite{le2015algorithm} for instance), but there is little theoretical justification for its use. Algorithm \ref{alg1} works by fixing some zeros of $F_{\delta}$ in Step 1, hence the relaxation (\ref{CNR}) while practically efficient, is not guaranteed to work, jeopardizing the entire construction.

\subsection{Critical Points Ideal}

A more modern relaxation comes from the gradient ideal $I_{\nabla} = \left \langle \frac{\partial p}{\partial x_{1}}, \dotso, \frac{\partial p}{\partial x_{n}} \right\rangle$. The first order optimality test $\nabla p(x) = 0$ implies that minima exist in the gradient variety $V_{\nabla}^{\mbb{R}}(I) =  \{ x\in \mbb{R}^{n} : \nabla p(x) = 0 \}$. In \cite{nie2006minimizing} it is shown that one may consider searching for minimizers in the quotient ring $\mbb{R}[x]/I_{\nabla}$ instead of $\mbb{R}[x]$. Their main theorem is the following;

\begin{theorem}[Theorem 8, \cite{nie2006minimizing}]\label{gradthm1}
	Assume that the gradient ideal $I_{\nabla}$ is radical. If the real polynomial $p(x)$ is non-negative over $V^{\mbb{R}}_{\nabla}(p)$, then there exist real polynomials $q_{i}(x)$ and $\phi_{j}(x)$ such that
	$$
	p(x) = \sum_{i=1}^{s}q_{i}(x) + \sum_{j=1}^{n} \phi_{j}(x)\frac{\partial p}{\partial x_{j}}(x)
	$$
	and each $q_{i}$ is a SOS.
\end{theorem}
Note that this is quite similar to (\ref{sosprog}), with the radicalness of $I_{\nabla}$ providing a guarantee on the existence of the decomposition. Algorithms for extracting the minimum and minimizers of functions are also presented in \cite{nie2006minimizing} and tested on several notable examples. In cases where it is unknown if $I_{\nabla}$ is radical, one may use the following alternative result of \cite{nie2006minimizing}.

\begin{theorem}[Theorem 9, \cite{nie2006minimizing}]\label{gradthm2}
	Suppose $p(x) \in \mbb{R}[x]$ is strictly positive on its real gradient variety $V_{\nabla}^{\mbb{R}}$. Then $p(x)$ is a SOS modulo its gradient ideal $I_{\nabla}$.
\end{theorem}

Extending Theorem \ref{gradthm1} and Theorem \ref{gradthm2}, \cite{demmel2007representations} considers the ideal generated by the KKT system related to $f$ when minimizing over a semialgebraic set. To this end let $\{ g_{1},\dotso,g_{j} \} \subseteq\mbb{R}[x]$ generate $S(G)$ and $PO(G)$. The KKT system associated to minimizing $p$ on $S(G)$ is 
\begin{equation*}
	\begin{gathered}
		\mc{P}_{i}  = \frac{\partial p}{\partial x_{i}} - \sum_{r=1}^{j} \lambda_{r} \frac{\partial g_{r}}{\partial x_{i}} = 0, \\
		g_{r}  \geq  0, \\
		\lambda_{r} g_{r}  \geq 0,
	\end{gathered}
\end{equation*}
for $r=1,\dotso,j$ and $i=1,\dotso,n$. As in \cite{demmel2007representations}, we let 
\begin{align*}
	I_{\text{KKT}} &= \left \langle \mc{P}_{1},\dotso,\mc{P}_{n}, \lambda_{1}g_{1}, \dotso, \lambda_{j}g_{j} \right\rangle, \\
	V^{\mbb{R}}_{\text{KKT}} & = \{ (x,\lambda) \in \mbb{R}^{n} \times \mbb{R}^{j} \ : \ q(x,\lambda)=0, \ \forall q\in I_{\text{KKT}}    \}, \\
	\mc{H} & = \{ (x,\lambda) \in \mbb{R}^{n} \times \mbb{R}^{j} \ : \ g_{r}(x,\lambda)\geq0, \ r=1,\dotso,j    \},
\end{align*}
and the KKT preorder generated by $G$ (now in the larger ring $\mbb{R}[x,\lambda]$) is 
$$
P_{\text{KKT}} = PO(G) + I_{\text{KKT}}.
$$ 

\begin{theorem}[Theorem 3.2, \cite{demmel2007representations}]
	Assume $I_{\text{KKT}}$ is radical. If $p(x)$ is non-negative on $V_{\text{KKT}}^{\mbb{R}} \cap \mc{H}$, then $p(x)$ belongs to $P_{\text{KKT}}$.
\end{theorem}
If the radicalness of $I_{\text{KKT}}$ is not known, then similar to Theorem \ref{gradthm1} positivity of $p(x)$ on the appropriate subset of $V_{\text{KKT}}^{\mbb{R}}$, ensures membership into $P_{\text{KKT}}$. 
\begin{theorem}[Theorem 3.5, \cite{demmel2007representations}]
	If $p(x)>0$ on $V_{\text{KKT}}^{\mbb{R}} \cap \mc{H}$, then $p(x)$ belongs to $P_{\text{KKT}}$.
\end{theorem}

For our application we work on the sphere $\mbb{S}^{N-1}$ (this can be replaced by any other suitable compact set) and the minimizers $(x^{*},y^{*})$ must now satisfy
\begin{equation*}
\label{kktsys}
	\begin{gathered}
	\nabla F_{\delta}(x,y) - \lambda \nabla s(x,y)  = 0, \\
	s(x,y)  = \sum_{i=1}^{n} x_{i}^{2} + \sum_{j=1}^{m} y_{j}^{2} - 1 = 0.
	\end{gathered}
\end{equation*}
This allows us to use the following KKT relaxation, 

\begin{equation}
\label{KKT}
	\begin{gathered}
		\max_{\delta>0}{\delta}  \\
		s.t. \ F_{\delta}(x) - \sum \phi_{i}(x) \left(\frac{\partial F_{\delta}}{\partial x_{i}}(x) -\lambda \frac{\partial s}{\partial x_{i}}(x)\right) -\lambda\eta(x)s(x) = SOS \\
		\phi_{i}, \eta \in \mbb{R}[x]
	\end{gathered}
\end{equation}

Notice that we do not search for membership of $F_{\delta}$ modulo $I_{\text{KKT}}$ into all of $PO(G)$, instead to simplify things we search only for elements of $PO(G)$ with $\gamma = (0,\dotsc,0)$. 
Since $F_{\delta}$ is known to have zeros, for this relaxation to be successful $I_{\text{KKT}}$ must be radical. While the random nature of $F_{\delta}$ implies a high probability of $I_{\text{KKT}}$ being radical, verifying this is computationally difficult, especially given the floating point construction of $F_{\delta}$. 


This relaxation also has non-linear constraints, arising from the products of decision variables (coefficients of $\phi_{i}$ and $\delta$). Hence, we implement this with the same ``bisection" approach and verification criteria as (\ref{HilB}). We fix $\delta = 2^{0}, k=1$, and solve
\begin{equation*}
\begin{gathered}
\text{find } \phi_{i}, \eta \in \mbb{R}_{k}[x,y], \\
s.t. \ F_{\delta}(x,y) - \phi(x,y)^{T}(\nabla F_{\delta}(x,y) - \lambda \nabla s(x,y)) -\lambda \eta(x,y) s(x,y) = SOS\\
\end{gathered}
\end{equation*}

\graphicspath{{Images/}}
\begin{figure}[H]
	\includegraphics[width=0.9\linewidth, center]{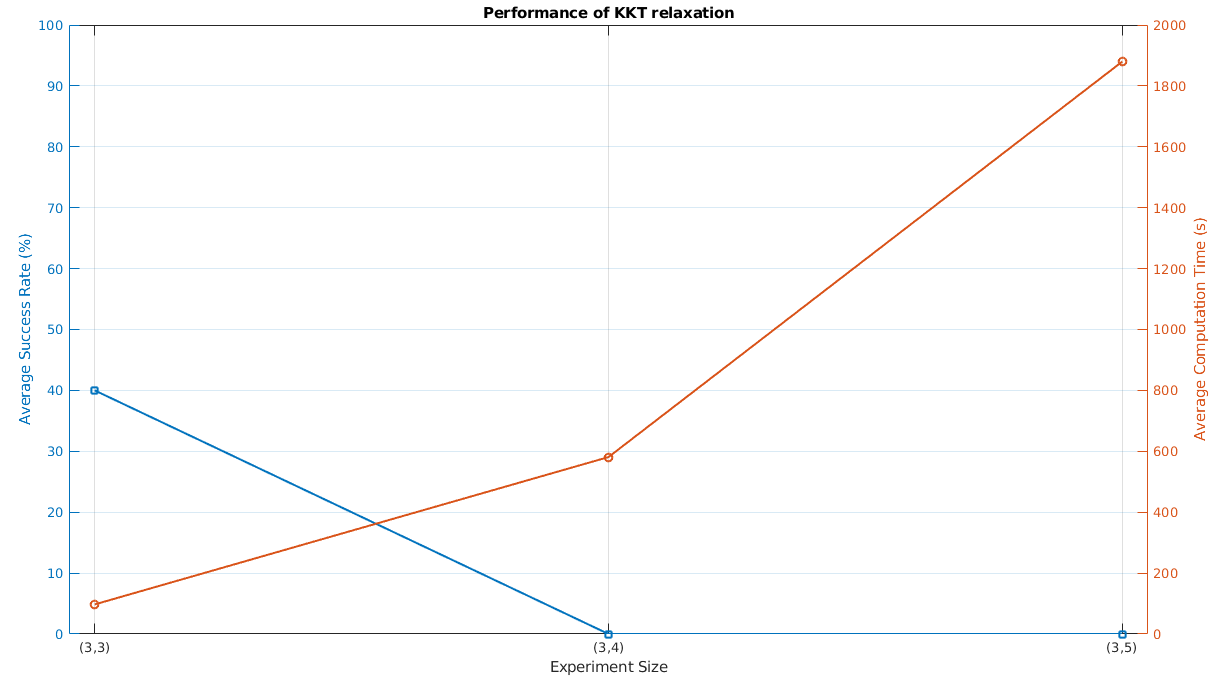}
	\caption{}
	\label{fig3}
\end{figure}
To our surprise, this method fails completely on the larger problems, and has quite poor performance even on the smaller ones of size $(3,3)$. This suggests that the random construction alone is not enough to guarantee the radicalness of $I_{\text{KKT}}$. Unlike the previous two relaxations, the residuals here do not indicate any room for improvement. In our tests, increasing the relaxation degree $k$ offers some success, but this also greatly increases the computation time, making this relaxation impractical for the problem at hand.

\begin{table}[H]
	\centering
	\begin{tabular}{ |p{2cm}||p{2cm}|p{2cm}|p{2cm}|  }
		\hline
		\multicolumn{4}{|c|}{KKT Relaxation}    \\
		\hline
		$(n,m)$ & Success (\%)  & Time (s)   & Residual   \\
		\hline
		$(3,3)$ & $40$          & $97.15$    & $16.95$     \\
		$(3,4)$ & $0$           & $581.06$   & $33.07$     \\
		$(3,5)$ & $0$           & $1879.94$  & $56.57$     \\
		\hline
	\end{tabular}
	\vspace{2pt}
	\caption{Average performance of relaxation (\ref{KKT})}
	\label{table:3}
\end{table}

\subsection{Jacobian relaxation}

We now present an exact relaxation which (in theory) always works for our problem of interest. This approach is similar to the KKT relaxation, only now to establish the dependence between derivatives of the constraints and the function, we consider determinants of an  associated Jacobian matrix.
Consider problems of the form (\ref{minprob}) with a single constraint $g$. Define the following
$$
		B(x) = \left[ \nabla p(x) \ \nabla g(x)  \right], \\
$$
and let
\begin{equation}\label{jacobian-phi}
		\varphi_{\ell}(x) = \sum_{\substack{E\in[N]_{2} \\ sum(E)=\ell}}\det B_{E}(x), \\
\end{equation}
where $B_{E}$ is the submatrix of $B$ with rows listed in $E$. 
As shown in \cite{nie2013exact}, (\ref{minprob}) is equivalent to 
\begin{equation}\label{jacobian}
	\begin{gathered}
		\min_{x\in\mbb{R}^{N}} p(x)\\
		s.t. \ \ g(x) = 0, \\
		\varphi_{\ell}(x) = 0, \ \ \ell = 3, \dotso, 2N-1.
	\end{gathered}
\end{equation}
We call this the \textit{Jacobian} system related to (\ref{minprob}). Letting $J = \langle g, \varphi_{3},\dotso,\varphi_{2N-1}\rangle$ and 
$$
J^{(k)} = \left\{ q\in J : deg(q)\leq2k \right\},
$$
we can write the SOS relaxation for Step 4 as
\begin{equation*}
	\begin{gathered}
		\max_{\gamma>0} \gamma, \\
		s.t. \ \ p(x) - \gamma = SOS + q(x), \\
		q(x) \in J^{(k)}.
	\end{gathered}
\end{equation*}

Moreover, letting $p^{*}$ be the solution of (\ref{jacobian}), $p^{(k)}$ of the corresponding SOS relaxation (of order $k$) and $p_{min}$ the minimum of (\ref{minprob}). Then the following holds.
\begin{theorem}[Theorem 2.3, \cite{nie2013exact}]
	Assume that $V(g)$ is non-singular, then $p^{*}>-\infty$ and there is a $K\in\mbb{N}$ such that $p^{(k)} = p^{*}$ for all $k\geq K$. Moreover, if $p_{min}$ is achievable, then $p^{(k)} = p_{min}$ for all $k\geq K$.
\end{theorem}
For us, the minimum of $F_{\delta}$ is always achieved on $\mbb{S}^{N-1}$, and it is clear that $V(s) = \mbb{S}^{N-1}$ is non singular. It follows that we can solve the Jocabian system (\ref{jacobian}) associated to $F_{\delta}$ exactly. This relaxation is given as 
\begin{equation}\label{JacRel}
	\begin{gathered}
		\max_{\delta>0} \delta, \\
		s.t. \ \ F_{\delta}(x) - q(x) = SOS, \\
		q(x) \in J^{(k)}.
	\end{gathered}
\end{equation}
Due to non-linearity in the constraints of (\ref{JacRel}), we employ the bisection approach similar to the other methods and solve
\begin{equation*}
	\begin{gathered}
		\text{find } q(x) \in J^{(k)}, \\
		s.t. \ \ F_{\delta}(x) - q(x) = SOS, 
	\end{gathered}
\end{equation*}
again with the limits of $\delta$ being $2^{-6}$ and $k$ being $2$.

\begin{remark}
	The functions $\varphi_{\ell}$ in (\ref{jacobian-phi}) are quartic polynomials in our problem of interest. The polynomials $q$ in (\ref{JacRel}) are also quartic polynomials. We could instead write this relaxation over the Segre variety in the variables $z_{ij} = x_{i}\otimes y_{j}$ which would lead to quadratic constraints $\varphi$. However, as detailed in \cite{nie2013exact} the generators of the Segre variety introduce an exponential number of constraints, and make (\ref{JacRel}) more difficult to solve numerically. This trade-off between the degree and the number of constraints is also present in the KKT relaxation. 
\end{remark}




\graphicspath{{Images/}}
\begin{figure}[H]
	\includegraphics[width=0.9\linewidth, center]{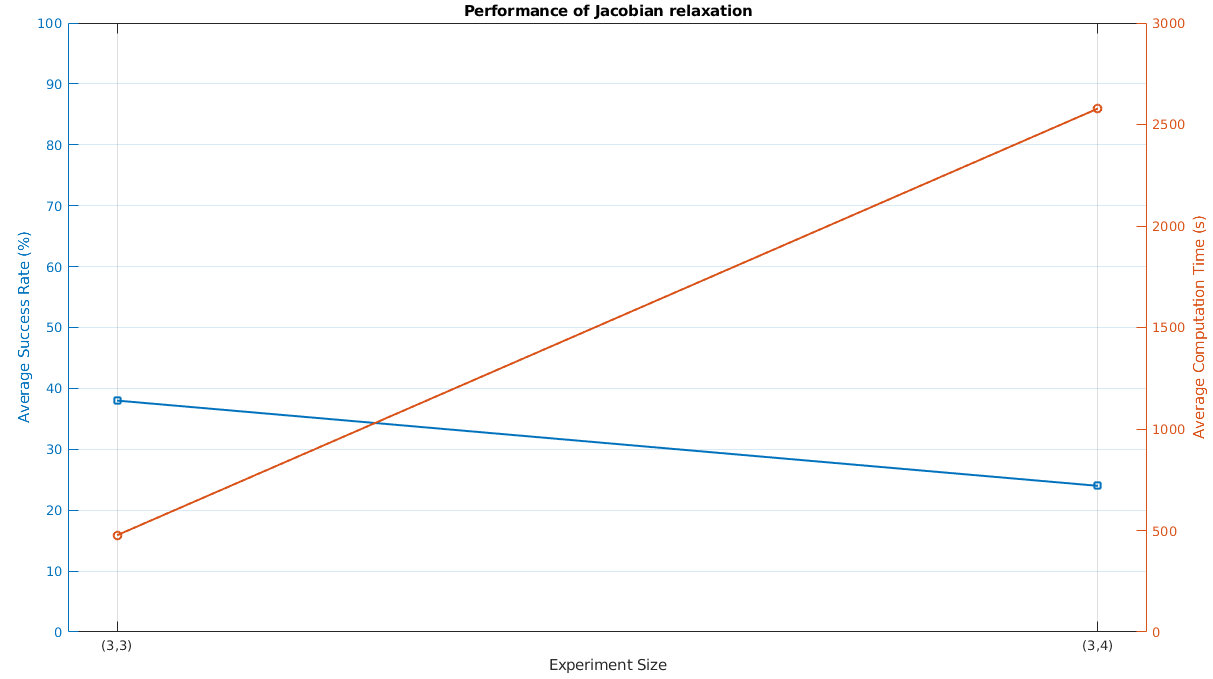}
	\caption{}
	\label{fig4}
\end{figure}

Unsurprisingly, this is quite slow. The solve time on test cases of size $(3,5)$ was close to one hour, and so we do not test the Jacobian relaxation on this set. We can also see (Figure/Table \ref{table:4}) that this relaxation exhibits low success rates and high residuals. Similar to KKT, the Jacobian relaxation is somewhat impractical in our context.

\begin{table}[H]
	\centering
	\begin{tabular}{ |p{2cm}||p{2cm}|p{2cm}|p{2cm}|  }
		\hline
		\multicolumn{4}{|c|}{Jacobian Relaxation}    \\
		\hline
		$(n,m)$ & Success (\%)    & Time (s)      & Residual    \\
		\hline
		$(3,3)$ & $38$            & $476.63$      & $18.64$     \\
		$(3,4)$ & $24$            & $2578.73$     & $25.06$     \\
		\hline
	\end{tabular}
	\vspace{2pt}
	\caption{Average performance of relaxation (\ref{JacRel})}
	\label{table:4}
\end{table}

\begin{remark}
	It should be noted again that these tests were conducted with limited freedom on the degrees of the relaxations. Based on our experience, we recommend using the Hilbert method with a high relaxation degree ($k=3$) if memory is not a concern and the user wants more successful constructions. When memory becomes an issue, the CNR cannot be beat; although its success rate is lower, the speed of computation makes generating random examples more practical.
\end{remark}

\section{Rationalization}

Constructing PnCP maps over floating point numbers provides quick numerical tests which can indicate non-negativity, but ideally we would like to have rational PnCP maps with exact certificates of non-negativity. The semi-definite programs arising from our SOS relaxations are feasibility problems of the form,
\begin{equation}\label{sdp}
	\begin{gathered}
		G\succeq 0 \\
		s.t. \quad \left\langle A_{i}, G \right\rangle  = b_{i}, \quad i=1,\dotso,m
	\end{gathered}
\end{equation}
where $A_{i}$ and $b_{i}$ are obtained from the problem data (see \cite{parrilo2003semidefinite} for a nice presentation of this). The following theorem, first proved in \cite{peyrl2008computing}, provides a means to obtain rational solutions of (\ref{sdp}) from numerical ones.

\begin{theorem}[Theorem 3.2, \cite{cafuta2015rational}]
	Let $G$ be a positive definite feasible point for (\ref{sdp}) satisfying
	$$
	\mu := \min{(eig(G))} > || (\left\langle A_{i}, G \right\rangle  - b_{i})_{i} || =: \epsilon, 
	$$
	then there is a (positive definite) rational feasible point $\hat{G}$. This can be obtained in two step;
	\begin{enumerate}[(1)]
		\item Compute a rational approximation $\tilde{G}$ with $\tau:=||G - \tilde{G}||$ satisfying $\tau^{2}+\epsilon^{2} \leq \mu^{2}$, 
		\item Project $\tilde{G}$ onto the affine subspace $\mc{L}$ defined by the equations $ \left\langle A_{i}, G \right\rangle = b_{i}$ to obtain $\hat{G}$. 
	\end{enumerate}
\end{theorem}

For our problems, there are two key issues with using this rationalization. Firstly, our SDP's will never satisfy the strict feasibility requirements of $G$ being positive definite. This is because by construction, the form $F_{\delta}$ will always have non-trivial zeros chosen in Step 1 of Algorithm \ref{alg1}. To tackle this, there are many facial reduction methods available to allow this rationalization for positive semi-definite $G$, one such reduction is presented in \cite{klep2017there} (see also \cite{laplagne2019facial} for instance).

More importantly, the numbers $b_{i}$ are obtained from the coefficients of the polynomial being tested, in our case $F_{\delta}$. This means that the affine subspace $\mc{L}$ is being defined by floating point numbers, and any sort of rationalization of $G$ will perturb this subspace. 

In PnCP we combat this by restricting the randomization in the linear algebra steps of Algorithm \ref{alg1}. As expected this reduces the base success rate of Algorithm \ref{alg1}, but it successfully constructs $F_{\delta}$ with rational coefficients. We also observe a significant increase in computation time to construct forms with rational coefficients; we test this by constructing 50 random forms with rational coefficients, and comparing the timing costs to constructing forms with floating point coefficients. 

As we can see below, constructing rational forms is far more expensive than floating point forms. In fact, the average time taken to construct forms with floating point coefficients remains almost constant ($\sim$2 seconds). In constrast, the construction time for forms with rational coefficients takes close to 10 minutes. 

\graphicspath{{Images/}}
\begin{figure}[H]
	\includegraphics[width=0.9\linewidth, center]{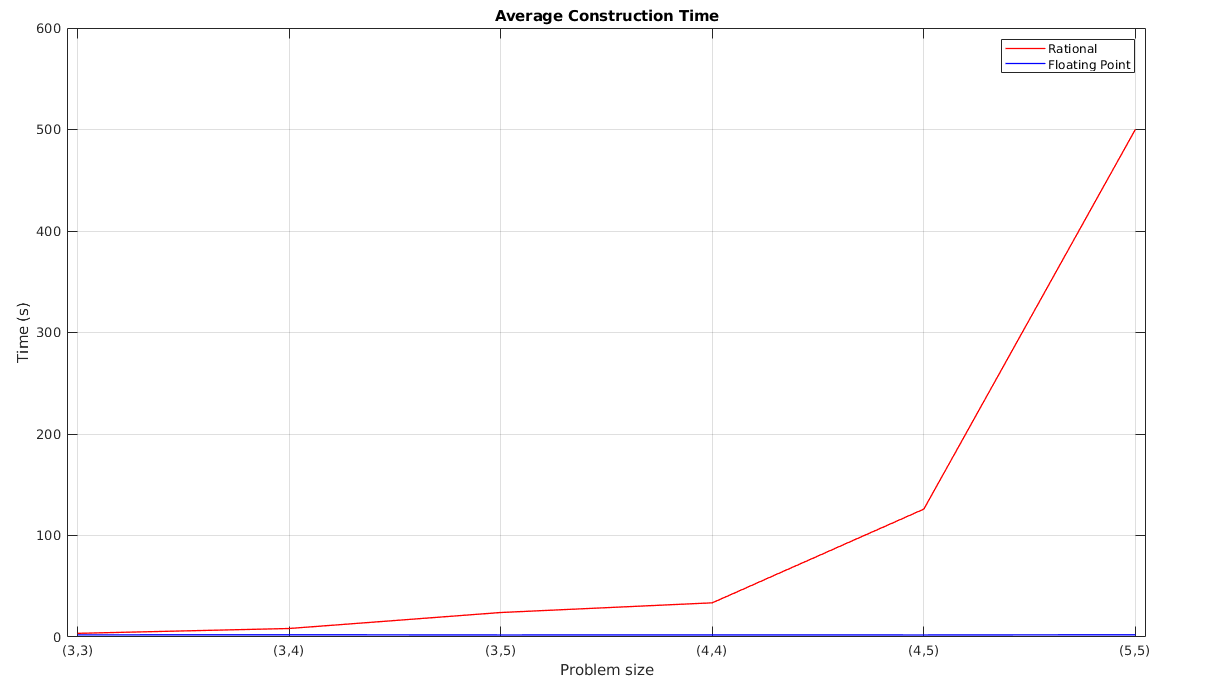}
	\caption{}
	\label{fig5}
\end{figure}

This rational construction can be used in PnCP with the command \texttt{Gen\_PnCP} and setting the \textit{`rationalize'} argument to 1. Currently, PnCP provides numerical verification of the constructed rational $F_{\delta}$, via the techniques of Section {\ref{relaxations}}. This construction can be used in conjunction with the many rational SOS packages (such as RationalSOS, RealCertify, multivsos, etc.) to obtain exact certificates of non-negativity.

\section{Detecting Quantum Entanglement}\label{quantum}

We will now show how we can use PnCP for detecting quantum entanglement. We start with a brief (and simplified) exposition into quantum states, the core object of interest for us, presenting some terminology and commonly known facts (for a more detailed introduction we refer the reader to \cite{vedral2006introduction, audretsch2008entangled, jaeger2007quantum}, or any graduate text on Quantum Information Theory). We then state two entanglement criteria, and then give an example demonstrating how PnCP is used to implement the most general one.

A quantum state is a vector $\phi \in \mbb{R}^{n}$, and with any quantum state there is an associated \textit{density matrix} $\phi\phi^{T}=:\rho\in M_{n}(\mbb{R})$. 
A density matrix
\begin{equation}\label{densitymat}
	\rho = \sum_{i} p_{i} \phi_{i} \phi_{i}^{T}, 
\end{equation}
with $\{\phi_{i}\}$ an orthonormal system, $p_{i}\geq0$ and $\sum_{i} p_{i} = 1$, represents a quantum system in one of several states $\phi_{i}$ with associated probabilities $p_{i}$. We use the following terminology; $\rho$ is a \textit{pure state} if $\rho = \phi \phi^{T}$, otherwise if $\rho$ is of the form (\ref{densitymat}), then it is a \textit{mixed state}. It should be noted that any positive semi-definite matrix $\rho$ with $Tr(\rho)=1$ is a density matrix. It is known that pure states satisfy $Tr(\rho^{2})=1$ while for mixed states $Tr(\rho^{2})<1$. 

Given a composite quantum system $M_{nm}(\mbb{R}) = M_{n}(\mbb{R})\otimes M_{m}(\mbb{R})$ and a state $\rho^{nm}\in M_{nm}(\mbb{R})$, we call $\rho^{nm}$ \textit{simply separable} if 
$$
\rho^{nm} = \rho^{n}\otimes\rho^{m}, \ \text{ with } \rho^{i}\in M_{i}(\mbb{R}), \text{ and } \norm{\rho^{i}}=1,
$$
\textit{separable} if
$$
\rho^{nm} = \sum_{i}p_{i}\rho_{i}^{n}\otimes\rho_{i}^{m}, \quad p_{i}\geq0, \quad \sum_{i} p_{i} = 1,
$$
and \textit{entangled} if its not separable. One of the big issues in quantum information theory is the so called \textbf{Separability Problem}; \textit{Given a state (density matrix) $\rho$ in a composite system, determine if it is entangled}.

There are many different criteria and measures of entanglement throughout the literature. For pure states, things are relatively simple and separability can be determined by checking if the state is in the image of the Segre embedding. For mixed states however, the situation is more complicated. 

In low dimensional composite systems, we have the Peres-Horodecki criterion, also known as the positive partial transpose (PPT) criterion; for
$
\rho^{nm} = \sum_{i}p_{i}\rho_{i}^{n}\otimes\rho_{i}^{m}
$
define the partial transpose map
$
(I\otimes T)(\rho^{nm}) = \sum_{i}p_{i}\rho_{i}^{n}\otimes(\rho_{i}^{m})^{T}.
$

\begin{criterion}[PPT, \cite{audretsch2008entangled} section 8.4]
	For a quantum state $\rho\in M_{nm}(\mbb{R})$, if $(I\otimes T)(\rho)$ has a negative eigenvalue, i.e., $(I\otimes T)(\rho)\nsucceq0$, then $\rho$ is entangled.
\end{criterion}

For systems of size $(n,m) = (2,2)$ or $(2,3)$, this criteria is both necessary and sufficient. In higher dimensional systems, we lose the sufficiency of this test, i.e., there are entangled states $\rho_{ent}$ with $(I\otimes T)(\rho_{ent})\succeq 0$ (see \cite{horodecki1997separability} for the first such example). In this situation we instead have the more general entanglement criteria.

\begin{criterion}[The general criterion, \cite{audretsch2008entangled} section 8.4]
	\label{entcrit}
	A quantum state $\rho\in M_{nm}(\mbb{R})$ is entangled if there is a pncp map $\Phi$ such that the ampliation $(I\otimes \Phi)(\rho)\nsucceq0$.
\end{criterion}
The PPT entanglement criterion is a special case of Criterion \ref{entcrit}, with $\Phi$ being the transposition map. With PnCP we can apply this test with many different random $\Phi$  in the following way;

\begin{algorithm}[H]
	\caption{Entanglement Detection}\label{alg2}
	\KwIn{$\rho$, $S$}
	\KwOut{Status}
	$i=0$\;
	Status = ``Unknown"\;
	\While{$i<S$}{
		Generate random $\Phi$\;
		Compute $I\otimes\Phi(\rho)$\;
		\eIf{$I\otimes\Phi(\rho)\nsucceq0$}{
			Status = ``Entangled"\;
			break\;
		}{
			$i = i+1$\;
		}
	}
\end{algorithm}

\begin{example}\label{example1}
	As an example consider the following state,
\begin{align*}
	\varDelta & = \begin{bmatrix}
		 1/3 & 0 & 0 & 0 & 1/3 & 0 & 0 & 0 & 1/3 \\
		   0 & 0 & 0 & 0 &   0 & 0 & 0 & 0 &   0 \\
		   0 & 0 & 0 & 0 &   0 & 0 & 0 & 0 &   0 \\
		   0 & 0 & 0 & 0 &   0 & 0 & 0 & 0 &   0 \\
		 1/3 & 0 & 0 & 0 & 1/3 & 0 & 0 & 0 & 1/3 \\
		   0 & 0 & 0 & 0 &   0 & 0 & 0 & 0 &   0 \\
		   0 & 0 & 0 & 0 &   0 & 0 & 0 & 0 &   0 \\
		   0 & 0 & 0 & 0 &   0 & 0 & 0 & 0 &   0 \\
		 1/3 & 0 & 0 & 0 & 1/3 & 0 & 0 & 0 & 1/3 
	\end{bmatrix} \\
	& = \frac{1}{3}\sum_{i,j} E_{i,j}\otimes E_{i,j} \in M_{3}(\mbb{R})\otimes M_{3}(\mbb{R}) ,
\end{align*}
where each $E_{i,j}\in M_{3}(\mbb{R})$ is the matrix unit with $1$ in row $i$, column $j$ and zeros everywhere else. This state is modeled after the Bell states, and is entangled. We use PnCP to generate the following non-negative, non-SOS polynomial with the command \texttt{Ent\_PnCP},
\begin{align*}
	F_{1}(x,y) & = 5x_{1}^{2}y_{1}^{2}+4x_{1}^{2}y_{1}y_{3}+12x_{1}x_{2}y_{1}^{2}-22x_{1}x_{2}y_{1}y_{2}+36x_{1}x_{2}y_{1}y_{3}+8x_{1}x_{3}y_{1}^{2}\\
	&\mathrel{\phantom{=}}+2x_{1}x_{3}y_{1}y_{2}+6x_{1}x_{3}y_{1}y_{3}+2x_{1}^{2}y_{2}^{2}+2x_{1}^{2}y_{2}y_{3}+60x_{1}x_{2}y_{2}^{2}-74x_{1}x_{2}y_{2}y_{3}\\
	&\mathrel{\phantom{=}}+4x_{1}x_{3}y_{2}^{2}+2x_{1}x_{3}y_{2}y_{3}-3x_{1}^{2}y_{3}^{2}+28x_{1}x_{2}y_{3}^{2}-2x_{1}x_{3}y_{3}^{2}+19x_{2}^{2}y_{1}^{2}\\
	&\mathrel{\phantom{=}}-66x_{2}^{2}y_{1}y_{2}+24x_{2}^{2}y_{1}y_{3}-4x_{2}x_{3}y_{1}^{2}+24x_{2}x_{3}y_{1}y_{2}-10x_{2}x_{3}y_{1}y_{3}+94x_{2}^{2}y_{2}^{2}\\
	&\mathrel{\phantom{=}}-36x_{2}^{2}y_{2}y_{3}+30x_{2}x_{3}y_{2}^{2}+2x_{2}x_{3}y_{2}y_{3}+5x_{2}^{2}y_{3}^{2}-2x_{2}x_{3}y_{3}^{2}+3x_{3}^{2}y_{1}^{2}\\
	&\mathrel{\phantom{=}}+2x_{3}^{2}y_{1}y_{2}+2x_{3}^{2}y_{1}y_{3}+2x_{3}^{2}y_{2}^{2}+x_{3}^{2}y_{3}^{2}
\end{align*}
and the associated PnCP map $\Phi$,

\begin{align*}
	\Phi(E_{1,1}) = \begin{bmatrix}
	     5  &   0  &   2 \\
		 0  &   2  &   1 \\
		 2  &   1  &  -3 
	\end{bmatrix}, & \quad
		\Phi(E_{1,3} +E_{3,1}) = \begin{bmatrix}
	8 & 1 &  3 \\
	1 & 4 &  1 \\
	3 & 1 & -2
	\end{bmatrix} \\
	\Phi(E_{3,3}) = \begin{bmatrix}
		3 & 1 & 1 \\
		1 & 2 & 0 \\
		1 & 0 & 1
	\end{bmatrix}, & \quad
	\Phi(E_{1,2} +E_{2,1}) = \begin{bmatrix}
	  12 & -11 &  18 \\
	 -11 &  60 & -37 \\
	  18 & -37 &  28
	\end{bmatrix}\\
	\Phi(E_{2,2}) = \begin{bmatrix}
	19 & -33 &  12 \\
	-33 &  94 & -18 \\
	12 & -18 &   5
	\end{bmatrix}, & \quad
	\Phi(E_{2,3}+E_{3,2}) = \begin{bmatrix}
		 -4 & 12 & -5\\
		 12 & 30 &  1\\
		 -5 &  1 & -2
	\end{bmatrix}
\end{align*}

Since we construct $\Phi$ on $\mbb{S}_{3}$, we make the canonical extension to $M_{3}(\mbb{R})$ by setting $\Phi(E_{i,j}) = \frac{1}{2}\Phi(E_{i,j}+E_{j,i})$ for $i\neq j$. With this extension, we find that
$$
(I\otimes\Phi)(\varDelta) = \frac{1}{6} \begin{bmatrix}
  10 &    0 &    4 &   12 &  -11 &   18 &   8 &   1 &   3 \\
   0 &    4 &    2 &  -11 &   60 &  -37 &   1 &   4 &   1 \\
   4 &    2 &   -6 &   18 &  -37 &   28 &   3 &   1 &  -2 \\
  12 &  -11 &   18 &   38 &  -66 &   24 &  -4 &  12 &  -5 \\
 -11 &   60 &  -37 &  -66 &  188 &  -36 &  12 &  30 &   1 \\
  18 &  -37 &   28 &   24 &  -36 &   10 &  -5 &   1 &  -2 \\
   8 &    1 &    3 &   -4 &   12 &   -5 &   6 &   2 &   2 \\
   1 &    4 &    1 &   12 &   30 &    1 &   2 &   4 &   0 \\
   3 &    1 &   -2 &   -5 &    1 &   -2 &   2 &   0 &   2
\end{bmatrix},
$$
with eigenvalues $  -8.45,   -2.78,   -0.83,   -0.06,    0.23,    2.17,    3.28,    7.14,   41.96$. 
\end{example}

\begin{example}\label{second-example}
	We consider now an example of a \textit{Bound Entangled State}, which are known to be entangled whilst having a positive partial transpose (see \cite{horodecki1998mixed} or \cite[Section 6.11]{jaeger2007quantum}). We take the example from \cite{halder2019construction}, with 
	\begin{equation*}
		\begin{gathered}
		\sigma  = \frac{1}{60}\begin{bmatrix}
						  5 &   5 &  -1 &  -1 &  -1 &  -1 &  -1 &  -1 &  -1 &  -1 &  -1 &  -1 \\
						  5 &   5 &  -1 &  -1 &  -1 &  -1 &  -1 &  -1 &  -1 &  -1 &  -1 &  -1 \\
						 -1 &  -1 &   5 &  -1 &  -1 &  -1 &  -1 &  -1 &  -1 &  -1 &  -1 &   5 \\
						 -1 &  -1 &  -1 &   5 &  -1 &  -1 &  -1 &  -1 &  -1 &   5 &  -1 &  -1 \\
						 -1 &  -1 &  -1 &  -1 &   5 &   5 &  -1 &  -1 &  -1 &  -1 &  -1 &  -1 \\
						 -1 &  -1 &  -1 &  -1 &   5 &   5 &  -1 &  -1 &  -1 &  -1 &  -1 &  -1 \\
						 -1 &  -1 &  -1 &  -1 &  -1 &  -1 &   5 &  -1 &   5 &  -1 &  -1 &  -1 \\
						 -1 &  -1 &  -1 &  -1 &  -1 &  -1 &  -1 &   5 &  -1 &  -1 &   5 &  -1 \\
						 -1 &  -1 &  -1 &  -1 &  -1 &  -1 &   5 &  -1 &   5 &  -1 &  -1 &  -1 \\
						 -1 &  -1 &  -1 &   5 &  -1 &  -1 &  -1 &  -1 &  -1 &   5 &  -1 &  -1 \\
						 -1 &  -1 &  -1 &  -1 &  -1 &  -1 &  -1 &   5 &  -1 &  -1 &   5 &  -1 \\
						 -1 &  -1 &   5 &  -1 &  -1 &  -1 &  -1 &  -1 &  -1 &  -1 &  -1 &   5
				 \end{bmatrix}, \\
			  \in M_{4}(\mbb{R})\otimes M_{3}(\mbb{R}).
		\end{gathered}
	\end{equation*}
	Note that $Tr(\sigma^{2}) = 0.2<1$, and so $\sigma$ is a mixed state (meaning we cannot simply check if it is in the image of the Segre embedding). PnCP generates the following
	
	\begin{align*}
	\Phi(E_{1,1}) = \begin{bmatrix}
		    7 & 17/2 & -5/2 \\
		 17/2 & 13/2 & -7/2 \\
		 -5/2 & -7/2 &   2 
	\end{bmatrix}, & \quad
	\Phi(E_{1,3} +E_{3,1}) = \begin{bmatrix}
		 -6 & -3 & 3 \\
		 -3 & -2 & 3 \\
		  3 &  3 & 0
	\end{bmatrix}, \\
	\Phi(E_{3,3}) = \begin{bmatrix}
		  3 & -1 &  0 \\
		 -1 &  0 & -1 \\
		  0 & -1 &  3 
	\end{bmatrix}, & \quad
	\Phi(E_{1,2} +E_{2,1}) = \begin{bmatrix}
		 -1/2 &  15/2 &    -6 \\
		 15/2 &    15 & -17/2 \\
		   -6 & -17/2 &  9/2 
	\end{bmatrix},\\
	\Phi(E_{2,2}) = \begin{bmatrix}
		  3 &   0  & -1 \\
		  0 & 17/2 & -4 \\
		 -1 &  -4  & 3
	\end{bmatrix}, & \quad
	\Phi(E_{2,3}+E_{3,2}) = \begin{bmatrix}
		  2 & -3 &  0 \\
		 -3 & -2 &  3 \\
		  0 &  3 & -2
	\end{bmatrix}.
	\end{align*}
	We find the ampliation $(I\otimes \Phi)(\sigma)$ to be
	$$
	\frac{1}{120}
	\begin{bmatrix}
		  133 &   162 &  -101 &  -17 &  -18 &   13 &  -17 &  -18 &   13 &   19 &  -30 &   13 \\
		  162 &   308 &  -182 &  -18 &  -52 &   22 &  -18 &  -52 &   22 &  -30 &  -52 &   10 \\
		 -101 &  -182 &   129 &   13 &   22 &  -21 &   13 &   22 &  -21 &   13 &   10 &   15 \\
		  -17 &   -18 &    13 &  163 &   36 &  -29 &  -17 &  -18 &   13 &   67 &   84 &  -17 \\
		  -18 &   -52 &    22 &   36 &  104 &  -44 &  -18 &  -52 &   22 &   84 &   26 &  -20 \\
		   13 &    22 &   -21 &  -29 &  -44 &   51 &   13 &   22 &  -21 &  -17 &  -20 &    3 \\
		  -17 &   -18 &    13 &  -17 &  -18 &   13 &   67 &   36 &    7 &   19 &  -18 &    1 \\
		  -18 &   -52 &    22 &  -18 &  -52 &   22 &   36 &  104 &  -44 &  -18 &   50 &  -26 \\
		   13 &    22 &   -21 &   13 &   22 &  -21 &    7 &  -44 &   75 &    1 &  -26 &   15 \\
		   19 &   -30 &    13 &   67 &   84 &  -17 &   19 &  -18 &    1 &  139 &   72 &  -29 \\ 
		  -30 &   -52 &    10 &   84 &   26 &  -20 &  -18 &   50 &  -26 &   72 &  128 &  -80 \\
		   13 &    10 &    15 &  -17 &  -20 &    3 &    1 &  -26 &   15 &  -29 &  -80 &   75
	\end{bmatrix},
	$$
	with eigenvalues of $ -0.14 ,   0.00,    0.06,    0.10  ,  0.27 ,   0.37  ,  0.60 ,   0.79,    1.01,    1.81   , 2.76  ,  4.69$.
	For this example, PnCP took $\sim$10 seconds to numerically check the entanglement status of the state, with majority of the time spent constructing the rational $\Phi$. If we desired only an indication of entanglement, we could repeat this with $\Phi$ having floating point entries, and the whole process would be significantly quicker.
\end{example}
\begin{remark}
	With Example \ref{second-example}, PnCP only claims that the given state is entangled, it does not claim that $\sigma$ is bound entangled, i.e., it does not check whether $\sigma$ is distillable \cite{bennett1996purification}. Distillation of quantum states is beyond the scope of this article.
\end{remark}

There are many other entanglement criteria that rely on testing a condition with some PnCP map. As we can see from the examples,  PnCP provides a means to implement these criteria by being able to generate random (\textit{rational}) pncp maps.

\section{Conclusions \& Future Work}

In this article we present PnCP; a MATLAB package for constructing positive maps which are not completely positive, with a focus on the practicality of this construction and its application to testing entanglement of quantum states. 

PnCP is an open-source package available from \url{https://bitbucket.org/Abhishek-B/pncp/}. The package implements state of the art optimization techniques to numerically ensure positivity of the constructed maps. PnCP is even able to construct pncp maps with rational coefficients, which can be used in conjunction with existing software to obtain not only numerical, but exact certificates of positivity.

We use the KMSZ construction which additionally provides a priori knowledge of some of the zeros of the constructed polynomial. While there is work on optimizing polynomials with zeros \cite{powers2006quantitative,castle2011polya}, there are restrictions on the zeros in these methods. Whether it is possible to adapt the zeros of the KMSZ construction to suit these methods, is something we wish to study in the future. 

As the only package for this kind of construction, we intend to maintain and improve PnCP in various means; implementing better non-negativity tests as they become available, optimizing the existing code (perhaps even pursuing parallel computing where possible), and including more entanglement criteria to improve the classification of quantum states. 

Our primary focus moving forward will be to strengthen PnCP as a classification tool for quantum states; primarily by implementing a rational SOS decomposition method which will automatically provide exact certificates of positivity. 





\end{document}